\documentclass[11pt]{article}

\usepackage{graphicx}

\newtheorem {theo} {\bf Theorem} [section]
\newtheorem {prop} [theo] {\bf Proposition}

\newtheorem {lem} [theo] {\bf Lemma}
\newtheorem {defn} [theo] {\bf Definition}

\newtheorem {obs}  [theo]{\bf Remark}

\newtheorem {afir}  [theo]{\bf Claim}

\newcommand{\CaixaPreta}{\rule{2mm}{2mm}}
\newcommand{\qed}{\hfill\caixapreta \vspace{5mm}}

\newcommand{\pat}{\partial}
\newtheorem {REMCURSIVA} [theo] {\bf Remark}

\newcommand{\be}{\begin{eqnarray}}
\newcommand{\ee}{\end{eqnarray}}

\newcommand{\benn}{\begin{eqnarray*}}
\newcommand{\eenn}{\end{eqnarray*}}

\newcommand{\bse}{\begin{equation}}
\newcommand{\ese}{\end{equation}}

\newcommand{\bsenn}{\begin{displaymath}}
\newcommand{\esenn}{\end{displaymath}}

\newcommand{\C}{\mbox{I${\!\!\!}$C}}

\newcommand{\R}{\mbox{I${\!}$R}}
\newcommand{\Z}{\mbox{Z${\!\!}$Z}}

\newcommand{\coes}{\c c\~oes }

\newcommand{\rv}{{\bf \mbox {r}}}
\newcommand{\vv}{{\bf \mbox {v}}}
\newcommand{\uu}{{\bf \mbox {u}}}
\newcommand{\xx}{{\bf \mbox {x}}}

\newcommand{\pp}{{\bf \mbox {p}}}

\title{Continuous Symmetric Perturbations of Planar Power Law Forces}
 \author{C. Azev\^edo and P. Ontaneda\footnote{The first author was 
 partially supported by a CNPq doctorate grant. The second author was 
 partially supported by a research grant from CNPq, Brazil.} }
\date{}

\usepackage[all]{xy}

\begin{document}

\maketitle

\vspace{0,4cm}

\begin{abstract}
We show the existence of periodic solutions for continuous symmetric perturbations
of certain planar power law problems.     
    
\end{abstract}

\vspace{0,4cm}

In this paper we study  continuous symmetric perturbations of planar power law problems of the
form   $\stackrel{..} {\mbox{\rv}}\,=\, g({\mbox{\rv}}, \mu)$,  where
the unperturbed problem is
$\stackrel{..} {\mbox{\rv}} \,=g({\mbox{\rv}},0)=-\frac{\mbox{$\kappa$}}{|\!|\mbox{{\mbox{\rv}}} |\!|^{\alpha
+2}}{\mbox{\rv}}\, \,,\kappa> 0$, $0\leq \alpha.$
In particular, if $\alpha =1$, the unperturbed problem is Kepler's problem.\\

We prove the existence of periodic solutions of perturbed problems as above, close
to a given circular orbit of the unperturbed problem. We have two cases. When
$\alpha =1$ (that is, the unperturbed problem is Kepler's problem) we will require that
the perturbed problems are symmetric with respect to the $x$ and $y$ axes. For $\alpha \neq 1$,
we will require just one symmetry.\\

Here are the statements of our main results, for $\alpha =1$ and 
$\alpha\neq 1$:\\

\noindent(Notation: ${\mbox{\rv}}(t,{\mbox{\xx}},{\mbox{\vv}},\mu)$ denotes the
solution of $\stackrel{..} {\mbox{\rv}}\,=\, g({\mbox{\rv}}, \mu)$ with initial conditions
${\mbox{\rv}}(0,{\mbox{\xx}},{\mbox{\vv}},\mu)={\mbox{\xx}}$ and
$\dot{\mbox{\rv}}(0,{\mbox{\xx}},{\mbox{\vv}},\mu)={\mbox{\vv}}$.)\\

{\theo  Let  $C$ be a circle centered at the origin  $(0,0)$ of $\R^{2}$,
${\mbox{\xx}}_{0}\in C\cap \{ x-axis\}$ and let $U$ an
open neighborhood of $C$ of the form $C\subset U\subset (\R^{2}-\{ (0,0)\})$.\\
Let $a>0$ and $g : U\times (-a,a)\rightarrow\R^{2}$ continuous such 
that\\
(i) $g({\mbox{{\mbox{\rv}}}},0)=-\frac{\mbox{$\kappa$}}{|\!|{\mbox{{\mbox{\rv}}}}
|\!|^{3} }
{\mbox{{\mbox{\rv}}}},\,\,\,\kappa>0$,\\
(ii) $g({\mbox{\rv}},\mu )$ is $C^{1}$ in the variable ${\mbox{\rv}}\in
U$, for each $\mu ,$\\
(iii) for all $\mu $, $g$ is invariant (as a vector field)
by the reflections
$$\varphi_{1}(x,y)=(-x,y), \,\,\,\varphi_{2}(x,y)=(x,-y).$$
\quad  Then there is $\delta_{0}$, $0<\delta_{0}<a$, with the 
following property.  For each $\mu\in
(-\delta_{0},\delta_{0})$ there is a velocity ${\mbox{\vv}}_{\mu}$ 
such that the solution
${\mbox{\rv}}_{{}_{{\mathrm{\mbox{\vv}}}_{\mu}},{}_{{\mbox{$\mu$}}}}(t)
:=
{\mbox{\rv}}
(t,{\mbox{\xx}}_{0},{\mbox{\vv}}_{\mu},\mu)$
of
$\stackrel{..}{\mbox{{\mbox{\rv}}}} \,=g(\mbox{{\mbox{\rv}}},\mu)$
is periodic.
 Moreover,  given $\eta >0$ we can choose  $\delta_{0}>0$ such that

(1) the traces of these solutions
are simple closed curves, symmetric with respect to the x and y axes, 
and enclose the origin,

(2)  all velocities ${\mbox{\vv}}_{\mu}$ are $\eta$-close, $\mu < 
\delta_{0}$.\label{3.1.1}}\\

{\theo
 Let  $C$ be a circle centered at the origin  $(0,0)$ of $\R^{2}$,
${\mbox{\xx}}_{0}\in C\cap \{ x-axis\}$ and $U$ an
open neighborhood of $C$ of the form $C\subset U\subset (\R^{2}-\{ (0,0)\})$.\\
Let $a>0$ and $g : U\times (-a,a)\rightarrow\R^{2}$ continuous such that\\
(i) $g(\mbox{{\mbox{\rv}}},0)=\,-\frac{\mbox{$\kappa$}}
{|\!|\mbox{\rv}|\!|^{\alpha +2}}
{\mbox{\rv}},$ where
$\mbox{{\mbox{\rv}}}=(x,y),\,\,\kappa > 0,$ and $\,\,0\leq \alpha 
,\,\,\, \alpha\neq 1,$\\
(ii) $g({\mbox{\rv}},\mu )$ is $C^{1}$ in the variable ${\mbox{\rv}}\in
U$, for each $\mu ,$\\
(iii) for all $\mu $, $g$ is invariant (as a vector field)
by the reflection
$$\varphi(x,y)=(x,-y).$$
\quad  Then there is $\delta_{0}$, $0<\delta_{0}<a$, with the 
following property.  For each $\mu\in
(-\delta_{0},\delta_{0})$ there is a velocity ${\mbox{\vv}}_{\mu}$ 
such that the solution
${\mbox{\rv}}_{{}_{{\mathrm{\mbox{\vv}}}_{\mu}},{}_{{\mbox{$\mu$}}}}(t)
:= {\mbox{\rv}} (t,{\mbox{\xx}}_{0},{\mbox{\vv}}_{\mu},\mu)$
of $\stackrel{..}{\mbox{{\mbox{\rv}}}} \,=g(\mbox{{\mbox{\rv}}},\mu)$ is periodic.
Moreover, given $\eta >0$ we can choose  $\delta_{0}>0$ such that

(1) the traces of these solutions
are simple closed curves, symmetric with respect to the x-axis, 
and enclose the origin,

(2)  all velocities ${\mbox{\vv}}_{\mu}$ are $\eta$-close, $\mu < 
\delta_{0}$.\label{3.2.1}}\\

Note that in both theorems $g$ is required to be defined just in a neighborhood
of the circle $C$.\\

In general, we can not claim that  
${\mu}\mapsto{\mbox {\vv}}_{\mu}$ is continuous, but it is 
possible to prove that the map ${\mu}\mapsto\mbox{ \vv}_{\mu}$
satisfies a property a bit weaker than that of continuity. This is given in the 
addenda to the theorems above. These results were used in \cite{AO}.
For more details see sections 2 and 3.\\

For $\alpha =1$ a similar result is proved in \cite{V1} using Poincar\'e's analytic
continuation. Our result is more general in
the sense that we do not demand the existence of any first integral. Also, we require
$g$ to be just continuous. Still, the result in \cite{V1} requires a little less symmetry.\\

There are also some related results in the literature obtained using
Calculus of Variations. These methods work best for $\alpha \geq 2$ (strong force),
see \cite{AB}, \cite{B1}. There are some results for $0<\alpha <2$, $\,\alpha\neq 1,\,$ assuming more
hypotheses, but in all cases the form of the perturbations is more restrictive and
the domain must have a particular form, see \cite{AC1},  \cite{AC2}, \cite{B2}. There are also some
results for $\alpha =1$, which require less symmetry, see \cite{AC1}. All these results hold, in general,
for any dimension but demand more regularity and the perturbations are always potential.
 Note that our result holds also for non potential perturbations, where variational methods cannot 
be applied.\\

For $\alpha =0$, we could not find related results obtained using Calculus of Variations.
In \cite{V2}, the results of \cite{V1} are extended to $\alpha =0$, but in this
case we demand less symmetry, and, as before, we do not demand the existence of any first integral and require
$g$ to be just continuous.\\

It is interesting to note that the proofs of the theorems above, given in sections 2 and 3,
do not use neither Calculus of Variations nor analytic continuation, but simple and elementary geometric
arguments. The results of this paper were motivated by (and needed in) the study of the fixed
homogeneous circle problem \cite{AO}. \\

This paper has three sections and two appendices.
In the first section we give some preliminary results. The case $\alpha =1$ is treated in section
2, and the case $\alpha \neq 1$ is treated in section 3. In the appendices we present proofs of some
needed results for which we could not find an exact reference in the literature.\\

\section {Preliminaries.}

We begin presenting two lemmas that show that 
to find a periodic solution of a problem  with symmetries
it is enough to construct  only a piece of a solution, with certain 
properties.

{\lem Let  $\Omega\subset\R^{2}$ be an open set,  with  $\varphi\Omega=\Omega$,
where  $\varphi (x,y)=(x,-y).$ Let  $f:\Omega\rightarrow \R^{2}$, be invariant by
$\varphi$, that is, $f(\varphi {\mbox\pp})=\varphi f(\mbox\pp),$
$\,{\mbox\pp}\in\Omega.$ If
$\mbox{\rv}:[0,\tau]\rightarrow \Omega,\,$ $\tau>0,$
$\,\, {\mbox{\rv}}(t)=(x(t),y(t)),$  is a solution of
\begin{equation}
\stackrel{..}{\mbox \rv} \,=f(\mbox \rv)
\label{5.0}
\end{equation}

\noindent such  that  ${\mbox \rv}(0),\,{\mbox \rv}(\tau)\in\{ 
\,x\}$-axis,
$($that is, $y(0)=y(\tau)=0$$)$ and $\dot {\mbox \rv}(0), \dot {\mbox \rv}(\tau)$
are vertical $($that is, $\dot x(0)=\dot x(\tau)=0$$)$, then  the
extension $\bar{\mbox \rv}$ of ${\mbox \rv}$, defined by:

$$\begin{array}{l}
\bar {\mbox \rv}(t)\,\, =\,\,\left\{
\begin{array}{lll}
{\mbox \rv}(t-2n\tau),& t\in[2n\tau, (2n+1)\tau],& n\in\Z
\\&&\\
\varphi {\mbox \rv}(2n\tau -t), &t\in[(2n-1)\tau,2n\tau], &n\in\Z
\end{array}
\right.
\end{array}$$

\noindent  is a periodic solution of (\ref{5.0}) with period $2\tau$.  
Moreover, the trace of  $\bar {\mbox \rv}$ is symmetric  with respect  to the   $x$-axis.
\label{3.0.14}}

\vspace{0,35cm}
\noindent {\bf Proof:} Since ${\mbox \rv}(t),$ $\,t\in[0,\tau]$,  is a solution 
of (\ref{5.0}), we have that ${\mbox \rv}(t-a)\,$ and
$ {\mbox \rv}(a+\tau-t),
\,t\in[a,a+\tau]$, are  also  solutions  of (\ref{5.0}), for all  
$a\in\R$.
Since $f$  is invariant by $\varphi$, we have that $\varphi
{\mbox \rv}(a-t),\,t\in[a-\tau,a]$, is also a solution  of (\ref{5.0}).
Therefore, each part in the definition of $\bar {\mbox \rv}$  is a solution 
of (\ref{5.0}).\\

A direct calculation shows that these parts, and its first derivatives,
coincide at the endpoints of the intervals where they are defined. 
In this way, $\bar {\mbox \rv}$ is a well defined solution of (\ref{5.0}).
Moreover,  $\bar {\mbox \rv}(0)=\bar {\mbox \rv}
(2\tau)$ and $\dot{\bar
{\mbox \rv}}(0)=\dot{\bar {\mbox \rv}}(2\tau),$ and   the trace of 
$\bar {\mbox \rv}$ is
symmetric  with respect  to the $x$-axis (see figure below). \CaixaPreta

\vspace{0,15cm}

\begin{figure}[!htb]
\centering
\includegraphics[scale=0.4]{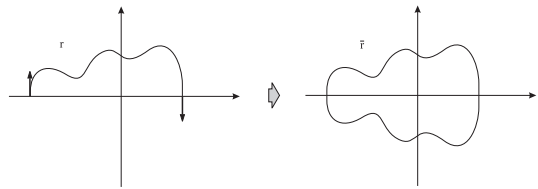}
\caption{\scriptsize{The solution is symmetric with respect to the $x$-axis}}
\end{figure}

\vspace{0,4cm}

Analogously we have:
{\lem Let  $\Omega\subset\R^{2}$ be an open set,  with 
$\varphi_{i}\Omega=\Omega$,
 $\,i=1,2$,
where   $\varphi_{1} (x,y)=(-x,y)$ and $\varphi_{2} (x,y)=(x,-y).$
Let  $f:\Omega\rightarrow \R^{2}$ be invariant by
$\varphi_{1}$ and $\varphi_{2}$. If
${\mbox \rv}:[0,\tau]\rightarrow \Omega,$ $\,\tau>0,\,$  is a solution of
\begin{equation}
\stackrel{..}{\mbox \rv}=f({\mbox \rv})
\label{5.a}
\end{equation}

\noindent such  that  ${\mbox \rv}(0)\in\{\,x\}$-axis, $\,{\mbox
\rv}(\tau)\in\{\,y\}$-axis,
$\dot {\mbox \rv}(0)$
is vertical and  $\dot {\mbox \rv}(\tau)$  is horizontal, then the
extension $\bar {\mbox \rv}$ of ${\mbox \rv}$, defined by:

$$\begin{array}{l}
\bar {\mbox \rv}(t)\,\, =\,\,\left\{
\begin{array}{lll}
{\mbox \rv}(t-4n\tau),& t\in[4n\tau, (4n+1)\tau],&
\\&&\\
\varphi_{1} {\mbox \rv}(4(n+2)\tau-t),
&t\in[(4n+1)\tau,(4n+2)\tau], &
\\&&\\
\varphi_{2}\varphi_{1} {\mbox \rv}(t-(4n+2)\tau),
&t\in[(4n+2)\tau,(4n+3)\tau],
&
\\&&\\
\varphi_{2} {\mbox \rv}(4n\tau-t), &t\in[(4n-1)\tau,(4n)\tau],
&
\end{array}
\right.
\end{array}$$

\noindent  with  $n\in\Z$,  is a periodic solution 
of (\ref{5.a}) with period $4\tau$. Moreover, the trace of  $\bar {\mbox \rv}$ is 
symmetric  with respect  to the $x$ and $y$ axes.
\label{3.0.15}}

\vspace{0,35cm}
\noindent {\bf Proof:} Analogous to the proof of the previous lemma (see figure 1.2). \CaixaPreta\\

\vspace{0,15cm}

\begin{figure}[!htb]
 \centering
\includegraphics[scale=0.4]{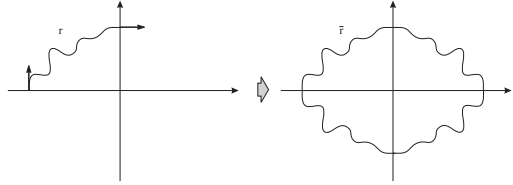}
\caption{\scriptsize{The solution is symmetric with respect to the  $x$ and $y$ axes}}
\end{figure}

\vspace{0,4cm}

We want to study now the transversality of solutions.

\vspace{0,1cm}

{\defn Let $I$ be an interval, $\alpha :I\rightarrow \R^{n}$
of class $C^{1}$ and $H^{n-1}$ a hypersurface of $\R^{n}$.
We say that  $\alpha$ intersects $H^{n-1}$ transversally if
$\alpha(\partial I)\cap H ^{n-1}=\emptyset,\,\,\,\alpha( I)\cap \partial H
^{n-1}=\emptyset,$ and $\,{\dot \alpha}(t)\notin T_{\alpha(t)}H^{n-1}$
for all   $t$  such  that  $\alpha(t)\in H^{n-1}$.
Moreover, we say that  $\alpha$ intersects 
$H^{n-1}$ transversally in a single point, if $\alpha$ intersects
$H^{n-1}$ transversally and there is a unique $t$ such that $\alpha(t)\in H^{n-1}.$}

\vspace{0,8cm}

The following proposition shows that  the property ``$\alpha$ 
intersects transversally in a single point''  is open in the  $C^{1}$ topology.

\vspace{0,25cm}

{\prop Let $E\subset\R^{2}$ be a closed segment   and let
$\alpha:[0,\bar t\,]\rightarrow\R^{2},\,$ $\bar t >0,$ be 
$C^{1}$, such that $\alpha$ intersects $E$ transversally in a single 
point.

Then there is  $\epsilon >0$ such that if $\beta :[0,\bar t\,]\rightarrow\R^{2}$
is $C^{1}$ and $\|\alpha -\beta\|<\epsilon,\,\,\|\dot\alpha
-\dot\beta\|<\epsilon,$ then   $\beta$ intersects $E$
transversally in a single point. (See figure 1.3.)
\label{3.0.17}}

\begin{figure}[!htb]
\centering
\includegraphics[scale=0.5]{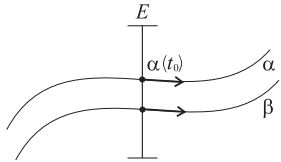}
\caption{\scriptsize{$E$ is closed, $\alpha$ and $\beta$ are transversal}}
\end{figure}

\vspace{0,3cm}

\noindent {\bf Remarks.}

\noindent (1) The fact that transversal maps form an open set in the $C^1$ topology
can be found in any differential topology textbook, but we could not
find a reference for the ``intersect in a single point" part of the
statement of the proposition above. Therefore we present a proof of this
proposition in an appendix.\\

\noindent (2) The  condition of $E$ being closed  is fundamental (see
figure below, where  $\alpha$ intersects transversally $E$ in a single
point  and $\beta$ is close to $\alpha,$ but $\beta$
intersects $E$ in two  points).

\begin{figure}[!htb]
 \centering
 \includegraphics[scale=0.5]{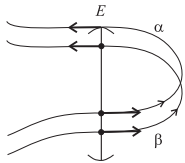}
\caption{\scriptsize{$E$ is not closed}}
 \end{figure}

\vspace{1cm}

The next proposition is essential in the proof of the theorems
in sections 2 and 3. Before, we introduce some notation:

\vspace{0,4cm}

Let   $U\times U_{\mu_{0}}\subset \R^{2}\times\R$ be an open set,  
where 
$\mu_{0}\in U_{\mu_{0}},$  and let
$g: U\times U_{\mu_{0}}\rightarrow \R^{2}$ be continuous and 
$C^{1}$ in $\mbox{{\mbox{\rv}}}\in U,$ for each $\mu$.
For each $\mu \in U_{\mu_{0}},$ consider the ordinary differential 
equation
\begin{equation}
 \stackrel{..}{\mbox{{\mbox{\rv}}}} \,= g({\mbox{{\mbox{\rv}}}},\mu )
\label{5.1}
\end{equation}

Denote by ${\mbox{\rv}}(t,{\mbox{\xx}},{\mbox{\vv}},\mu)$ a
solution of (\ref{5.1}), with initial conditions
${\mbox{\rv}}(0,{\mbox{\xx}},{\mbox{\vv}},\mu)={\mbox{\xx}}$ and
$\dot{\mbox{\rv}}(0,{\mbox{\xx}},{\mbox{\vv}},\mu)={\mbox{\vv}}$.
Let ${\mbox{\xx}}_{0}\in 
U,\,{\mbox{\vv}}_{0}\in\R^{2},\,\mu_{0}\in\R$ be fixed. 
We write $\mbox{{\mbox{\rv}}}  _{0}(t)=
{\mbox{\rv}}(t,{\mbox{\xx}}_{0},{\mbox{\vv}}_{0},\mu_{0}).$
In this situation we have:\\

 {\prop Let $[0,\bar t\,],$ $\,\bar t >0,$ be a
time interval in which $\mbox{{\mbox{\rv}}} _{0}( t)$
intersects transversally the closed segment  $E\subset \R^{2}$
in a single point.

 Then there is $\delta >0$ $\left( \delta \right.$ depending on
$({\mbox{\xx}}_{0},{\mbox{\bf{\mbox{\vv}}}}_{0},\mu_{0},\bar
t\,)\left.\right)$ such that for ${\mbox{\xx}},{\mbox{\vv}},\mu$
satisfying
$|\!|{\mbox{\xx}}-{\mbox{\xx}}_{0}|\!|<\delta,\,\,
|\!|{\mbox{\bf{\mbox{\vv}}}}-{\mbox{\bf{\mbox{\vv}}}}_{0}|\!|<\delta,\,\,
|\mu -\mu _{0}|<\delta,$
the solution  $\mbox{{\mbox{\rv}}}(t,{\mbox{\xx}},{\mbox{\bf{\mbox{\vv}}}},\mu)$
of the perturbed  problem (\ref{5.1})
is defined in $[0,\bar t\,]$ and, restricted to the interval $[0,\bar t\,],$
intersects  $E$ transversally in a single point.
Moreover, the function $0<t({\mbox{\xx}},{\mbox{\bf{\mbox{\vv}}}},\mu )<\bar t,$
defined by $\mbox{{\mbox{\rv}}} (t({\mbox{\xx}},{\mbox{\bf{\mbox{\vv}}}},\mu),
{\mbox{\xx}},{\mbox{\bf{\mbox{\vv}}}},\mu)\in E,$ is continuous.
\label{3.0.18}}

\vspace{0,4cm}

\noindent {\bf Proof:} 
Let $\epsilon>0$ be as in proposition  \ref{3.0.17} (taking $\alpha ={\mbox{\rv}}_{0}$).
Since solutions of ODE depend continuously on the initial data (e.g. see \cite{S}, p.34), there
exist $\delta>0\,$ ($\delta$
depending on $({\mbox{\xx}}_{0},{\mbox{\vv}}_{0},\mu_{0},\bar t\,))$ such that
if $({\mbox{\xx}},{\mbox{\vv}},\mu)$ satisfies:
\begin{equation}
 |\!|{\mbox{\xx}}-{\mbox{\xx}}_{0}|\!|<\delta,\,\,\,\,
|\!|{\mbox{\vv}}-{\mbox{\vv}}_{0}|\!|<\delta,\,\,\,\,|\mu-\mu_{0}|<\delta,
\label{5.2}
\end{equation}
then the solution  $\mbox{{\mbox{\rv}}}=
\mbox{{\mbox{\rv}}}(t,{\mbox{\xx}},{\mbox{\mbox{\vv}}},\mu)$ of the
problem  (\ref{5.1}) is defined on $[0,\bar t\,]$ and satisfies: 
\vspace{0,15cm}
$$\begin{array}{lll}
\|\mbox{{\mbox{\rv}}}-\mbox{{\mbox{\rv}}}_{0}\|&=& sup_{t\in I}
\|{\mbox{\rv}}(t,{\mbox{\xx}},{\mbox{\vv}},\mu)-{\mbox{\rv}}_{0}
(t,{\mbox{\xx}}_{0},{\mbox{\vv}}_{0},\mu_{0})\|<\epsilon;
\\&&\\
\|\dot{\mbox{\rv}}-{\dot{\mbox{\rv}}}_{0}\|&=& sup_{t\in I}
\|
 \dot{\mbox{\rv}}(t,{\mbox{\xx}},{\mbox{\vv}},\mu)-
{\dot{\mbox{\rv}}} _{0}
(t,{\mbox{\xx}}_{0},{\mbox{\vv}}_{0},\mu_{0})\|< \epsilon.
\end{array}$$

By proposition  \ref{3.0.17}, we have  that   $\mbox{{\mbox{\rv}}}([0,\bar 
t\,])$ intersects $E$ transversally in a single point. Also, by this same proposition,
we have that the function $t({\mbox{\xx}},{\mbox{\vv}},\mu)$ is well defined for any
$({\mbox{\xx}},{\mbox{\vv}},\mu)$ satisfying (\ref{5.2}).
Write $({\mbox{\xx}},{\mbox{\vv}},\mu)={\mbox{\uu}}$.
Since the function $t$ is bounded, to prove its continuity it is enough to show that, 
if ${\mbox{\uu}}_{n}\rightarrow {\mbox{\uu}}$ and $t({\mbox{\uu}}_{n})\rightarrow
b$, then  $b=t({\mbox{\uu}}).$\\

By definition of $t$, we have  that  $\mbox{{\mbox{\rv}}} (t({\mbox{\uu}}_{n}),
{\mbox{\uu}}_{n})\in E.$ Since ${\mbox{\rv}}$ is continuous we have that 
$lim_{n\rightarrow +\infty}\mbox{{\mbox{\rv}}}(t({\mbox{\uu}}_{n}),{\mbox{\uu}}_{n})=
\mbox{{\mbox{\rv}}}(b,{ \mbox{\uu}})$; hence $\mbox{{\mbox{\rv}}}(b,{ \mbox{\uu}})\in E$
because $E$  is closed.
But $\mbox{{\mbox{\rv}}}(t({\mbox{\uu}}),{\mbox{\uu}})\in E$ and,
by proposition  \ref{3.0.17}, $t({\mbox{\uu}})$  is unique. This 
implies that $b=t({\mbox{\uu}}).$  This proves the proposition. \CaixaPreta

\vspace{0,5cm}

\section{Perturbations of Kepler's problem.}
 In this section we prove the existence of periodic solutions of perturbed problems 
with symmetries, close to a circular solution of Kepler's problem (the unperturbed problem).
\vspace{0,2cm}

As before, ${\mbox{\rv}}(t,{\mbox{\xx}},{\mbox{\vv}},\mu)$ denotes a
solution  of $\stackrel{..}{\mbox{{\mbox{\rv}}}} \,= g({\mbox{{\mbox{\rv}}}},\mu )$,
with initial conditions ${\mbox{\rv}}(0,{\mbox{\xx}},{\mbox{\vv}},\mu)={\mbox{\xx}}$ and
$\dot{\mbox{\rv}}(0,{\mbox{\xx}},{\mbox{\vv}},\mu)={\mbox{\vv}}$.\\

In this section we prove theorem 0.1:\\

\noindent {\bf Theorem 0.1} {\it  Let  $C$ be a circle centered at the origin  $(0,0)$ of $\R^{2}$,
${\mbox{\xx}}_{0}\in C\cap \{ x-axis\}$ and let $U$ an
open neighborhood of $C$ of the form $C\subset U\subset (\R^{2}-\{ (0,0)\})$.\\
Let $a>0$ and $g : U\times (-a,a)\rightarrow\R^{2}$ continuous such 
that\\
(i) $g({\mbox{{\mbox{\rv}}}},0)=-\frac{\mbox{$\kappa$}}{|\!|{\mbox{{\mbox{\rv}}}}
|\!|^{3} }
{\mbox{{\mbox{\rv}}}},\,\,\,\kappa>0,$\\
(ii) $g({\mbox{\rv}},\mu )$ is $C^{1}$ in the variable ${\mbox{\rv}}\in
U$, for each $\mu ,$\\
(iii) for all $\mu $, $g$ is invariant (as a vector field)
by the reflections
$$\varphi_{1}(x,y)=(-x,y), \,\,\,\varphi_{2}(x,y)=(x,-y).$$
\quad  Then there is $\delta_{0}$, $0<\delta_{0}<a$, with the 
following property.  For each $\mu\in
(-\delta_{0},\delta_{0})$ there is a velocity ${\mbox{\vv}}_{\mu}$ 
such that the solution
${\mbox{\rv}}_{{}_{{\mathrm{\mbox{\vv}}}_{\mu}},{}_{{\mbox{$\mu$}}}}(t)
:=
{\mbox{\rv}}
(t,{\mbox{\xx}}_{0},{\mbox{\vv}}_{\mu},\mu)$ of
$\stackrel{..}{\mbox{{\mbox{\rv}}}} \,=g(\mbox{{\mbox{\rv}}},\mu)$ is periodic.
Moreover,  given $\eta >0$ we can choose  $\delta_{0}>0$ such that

(1) the traces of these solutions
are simple closed curves, symmetric with respect to the x and y axes, 
and enclose the origin,

(2)  all velocities ${\mbox{\vv}}_{\mu}$ are $\eta$-close, $\mu < 
\delta_{0}$.}
\vspace{0,4cm}

\noindent {\bf Remarks 2.1}

\noindent (1) It can easily be deduced from proposition \ref{3.0.17}
and (2) of the theorem above that we can choose $\delta_{0}$ in the
theorem such that the periodic solutions intersect transversally the $x$-axis in 
exactly two points. Follows that these points are ${\mbox{\xx}}_{0}$
and $-{\mbox{\xx}}_{0}$. This fact is used in \cite{AO}.\\

\noindent (2)  Recall that a simple closed curve S, in the plane $\R^2,$ encloses
a point $p\notin S$, if $p$ belongs to the bounded component 
of $\R^2\setminus S$.\\

\noindent (3) $\delta_{0}$ of theorem \ref{3.1.1} depends only on
$g,\,$  $ U$ and the radius of $ C.$\\

\noindent (4) From the proof of the theorem (or from the fact that the solutions 
${\mbox{\rv}} _{ {} _{{\mbox{\vv}}_{\mu}},{}_{\mbox{$\mu$}}}$
are symmetric) follows  that  ${\mbox{\vv}}_{\mu}$  is vertical,
that is, ${\mbox{\vv}}_{\mu}=(0,\vv_{\mu})$ or, equivalently,
orthogonal to the $x-$axis.\\

\noindent (5) In general, we cannot claim that  
${\mu}\mapsto{\mbox {\vv}}_{\mu}$ is continuous, but it is 
possible to prove that the map ${\mu}\mapsto\mbox{ \vv}_{\mu}$
satisfies a property a bit weaker than that of continuity:
there is a correspondence $\mu \mapsto
{\cal V}_{\mu}\subset\R^{2},\,{\cal V}_{\mu}\neq \emptyset,$ such  that  $
{\mbox{\rv}}_{{}_{{\mbox\vv}},{}_{\mbox{$\mu$}}}$ is periodic, for all   $\mbox{ \vv}\in {\cal V}
_{\mu}.$  Moreover, the set ${\cal V}=\cup_{\mu\in[0,\delta_{0}]}{\cal V}_{\mu}=\{
({\mbox{\vv}},{\mu});\,\mbox{\vv}\in {\cal V}_{\mu}\}$  is compact and connected. (See the figure below.)

\begin{figure}[!htb]
\centering
\includegraphics[scale=0.5]{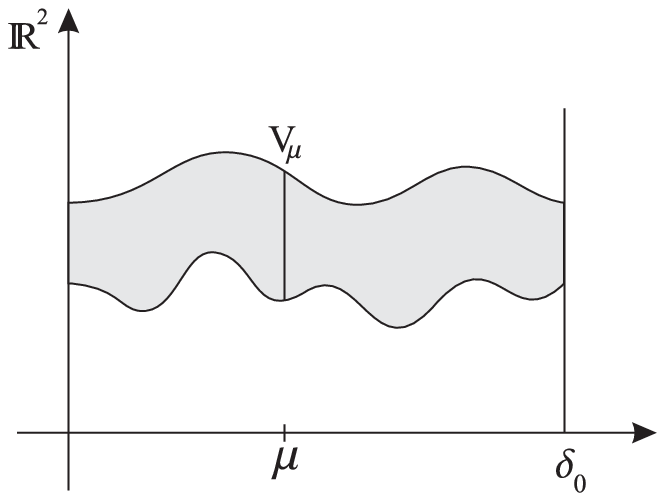}
\end{figure}

\vspace{0,2cm}

If  ${\mu}\mapsto{\mbox {\vv}}_{\mu}$
is already continuous, we can take ${\cal V}_{\mu}=\{ {\mbox {\vv}}_{\mu}\}.$ In this particular case 
${\cal V}=\cup_{\mu\in[0,\delta_{0}]}{\cal V}_{\mu}=\{
({\mbox{\vv}},{\mu});\,\mbox{\vv}\in {\cal V}_{\mu}\}$ 
is the graph of 
${\mu}\mapsto{\mbox {\vv}}_{\mu}$, $\mu\in [0,\delta_{0}]$, which is certainly compact and connected.\\

We make then the following addendum to theorem \ref{3.1.1}.\\

\noindent {\bf Addendum (to theorem \ref{3.1.1})} {\it 
We can choose $\delta_{0}>0$ in theorem \ref{3.1.1} such that there is a compact connected
${\cal V}\subset \R^{2}\times\R$  with the following properties:

\noindent 1) ${\cal V}_{\mu}:={\cal V}\cap \left( \R^{2}\times\{\mu\}\right)\neq
\emptyset,$ for all   $\mu\in[0,\delta_{0}]$,

\noindent 2) ${\mbox{\rv}}_{{}_{\mbox{\vv}},{}_{\mbox{$\mu$}}}$ is periodic, for
$(\mbox{\vv},\mu)\in \cal V.$

Moreover, the trace of ${\mbox{\rv}}_{{}_{\mbox{\vv}},{}_{\mbox{$\mu$}}}$ is a simple closed curve 
symmetric with respect to the x and y axes, and encloses the origin.}

\vspace{0,5cm}

\noindent {\bf Proof of the Theorem \ref{3.1.1}:}
First, recall that every circle centered at the origin is the trace of  a 
periodic solution of Kepler's problem. This solution  has constant angular speed $w$,
with  $w=\sqrt{\kappa}a^{-3/2},$ where $a$  is the radius of the
circular solution.\\

Let   ${\mbox{{\mbox{\rv}}}}_{0}(t)$ be the circular solution of Kepler's problem in
the $(x,y)$-plane whose trace is $C$. We can assume that ${\mbox{{\mbox{\rv}}}}_{0}(0)={\mbox{\xx}}_{0}$
lies in the positive $x$-axis and  $\dot{{\mbox{{\mbox{\rv}}}}}_
{0}(0)={\mbox{\vv}}_{0}\,$ has the same direction as the positive
$y$-axis, that is, ${\mbox{\vv}}_{0}=(0,{{\vv}}_{0}),\, {{\vv}}_{0}>0$.\\

Choose $\bar t >0$ such  that  ${\mbox{\rv}}_{0}(t)=
{\mbox{\rv}}(t,{\mbox{\xx}}_{0},{\mbox{\vv}}_{0},0)$,
restricted to the interval $[0,\bar t\,]$, intersects the
$y$-axis transversally in a single point (see figure below).

\begin{figure}[!htb]
\centering
\includegraphics[scale=0.45]{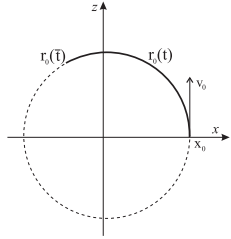}
\caption{\scriptsize{${\mbox{\rv}}_{0}(t),\,\,t\in [0,\bar t\,],$
intersects the $y$-axis  transversally  in a single   point }}
\end{figure}

\vspace{0,4cm}

Consider $\delta$, $0<\delta<1$, $\delta$
sufficiently small (as in proposition  \ref{3.0.18},  with  $\bar t$  as above and
$E=$$y$-axis), and define
$V_{\delta} := \{\,\sigma {\mbox{\vv}}_{0}\,; \,\sigma\in(1-\delta,
1+\delta)\,\}\subset\R^{2}$.
Let $l:V_{\delta}\times (-\delta,\delta)\rightarrow
\R,\,\,l({\mbox{{\mbox{\vv}}}},\mu)=\,\dot{y}(t({\mbox{{\mbox{\vv}}}},\mu),{\mbox{\xx}}_{0},
{\mbox{\vv}},\mu),$ where 
$\dot{y}(t({\mbox{{\mbox{\vv}}}},\mu),{\mbox{\xx}}_{0},{\mbox{\vv}},\mu)$  is the second coordinate of
$\dot{\mbox{{\mbox{\rv}}}}(t({\mbox{{\mbox{\vv}}}},\mu),{\mbox{\xx}}_{0},{\mbox{\vv}},\mu)
=(\dot{{\mbox{\xx}}}(t({\mbox{{\mbox{\vv}}}},\mu),
{\mbox{\xx}}_{0},{\mbox{\vv}},\mu),\dot{y}(t({\mbox{{\mbox{\vv}}}},\mu),{\mbox{\xx}}_{0},{\mbox{\vv}},\mu))$
and $t({\mbox{\vv}},\mu)$
is the time at which the  solution  intersects the $y$-axis.
Here $t({\mbox{\vv}},\mu):= t({\mbox\xx}_{0},{\mbox{\vv}},\mu)$  is as in
proposition  \ref{3.0.18} (this is why we chose $\bar t$ and wanted
$\delta$ sufficiently small). Note that, since $\delta<1$,
${\mbox{\vv}}\in V_{\delta}$ has the same direction as ${\mbox{\vv}}_{0}$.\\

Observe that   $l({\mbox{{\mbox{\vv}}}},\mu)=0$ if and only if
${\mbox{\rv}}_{{}_{\mbox{\vv}},{}_{\mbox{$\mu$}}}$,
restricted to $[0,\bar t\,]$, intersects orthogonally the $y$-axis
(and in a single point).\\

{\afir {\rm  Let  $\eta>0$. Then there is 
$\delta_{0},$ $0<\delta _{0}<\delta$,
and ${\mbox{\vv}}_{+}, {\mbox{\vv}}_{-}\in V_{\delta}$
such that 
$|\!| {\mbox{\vv}}_{0}-{\mbox{\vv}}_{-}|\!|<\eta$,
$|\!| {\mbox{\vv}}_{0}-{\mbox{\vv}}_{+}|\!|<\eta$,
$|\!| {\mbox{\vv}}_{-}|\!|<|\!|{\mbox{\vv}}_{0}|\!|<|\!|{\mbox{\vv}}_{+}|\!|$,
and for all   $\mu\in (-\delta_{0},\delta_{0}),$
we have $l({\mbox{\mbox{\vv}}}_{+},\mu)>0$
and $l({\mbox{\mbox{\vv}}}_{-},\mu)<0$.
\label{af}}}

\vspace{0,3cm}

\noindent {\bf Proof of the claim:}
For $\mu =0$ we have Kepler's problem. Hence there is ${\mbox{\vv}}_{+}$ and
${\mbox{\vv}}_{-}$ belonging to $V_{\delta}$, such that 
$|\!| {\mbox{\vv}}_{-}|\!|<|\!|{\mbox{\vv}}_{0}|\!|<|\!|{\mbox{\vv}}_{+}|\!|$ and
$|\!| {\mbox{\vv}}_{0}-{\mbox{\vv}}_{-}|\!|<\eta$,
$|\!| {\mbox{\vv}}_{0}-{\mbox{\vv}}_{+}|\!|<\eta$,  with 
$l({\mbox{\vv}}_{+},0)>0$ and
$l({\mbox{\vv}}_{-},0)<0$ (see figure below).
\vspace{0,2cm}

\begin{figure}[!htb]
\centering
\includegraphics[scale=0.45]{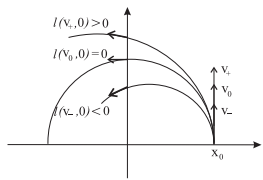}
\caption{\scriptsize{Solutions close to C}} 
\end{figure}

\vspace{0,3cm}

Since $l$  is continuous, for $\epsilon
=\frac{1}{2}\,min\{l({\mbox{\vv}}_{+},0),\, |l({\mbox{\vv}}_{-},0)|\}>0$, there is
$\delta_{0}>0$ so  that  $0<\delta_{0}<\delta$ and
$|l({\mbox{\vv}}_{-},0)-l({\mbox{\vv}}_{-},\mu)|<\epsilon$,
$|l({\mbox{\vv}}_{+},0) -l({\mbox{\vv}}_{+},\mu)|<\epsilon$, for
$\mu\in(-\delta_{0},\delta_{0})$.
It follows that 
$l({\mbox{\vv}}_{+},\mu)>0,$ and $l({\mbox{\vv}}_{-},\mu)<0,$  for
$\mu\in(-\delta_{0},\delta_{0})$. This proves the claim.

\vspace{0,6cm}

By the claim above (choose any $\eta >0$) and by the intermediate value theorem, there is
$\delta_{0}>0$ such  that for each $\mu\in(-\delta_{0},\delta_{0})$
there exists ${{\mbox{\vv}}}_{\mu}\in V_{\delta}$ satisfying
$l({{\mbox{\vv}}}_{\mu},\mu)=0$. Therefore, for each  $\mu,$ the solution of the
perturbed problem 
${\mbox{\rv}}_{\mu}(t):={\mbox{\rv}}_{{}_{{\mbox\vv}_{\mu}},{}_{\mbox{$\mu$}}}(t)
=
{\mbox{\rv}}(t,{\mbox{\xx}}_{0},{\mbox{\vv}}_{\mu},\mu)$
intersects orthogonally the $y$-axis at $t=t({{\mbox{\vv}}}_{\mu},\mu).$\\

Finally, lemma \ref{3.0.15} implies that 
${\mbox{\rv}}_{\mu}(t)$ can be extended to a periodic solution, with  period $4t({{\mbox{\vv}}}_{\mu},\mu)$, 
whose trace  is symmetric  with respect  to the $x$ and $y$ axes.\\

To show that we can choose $\delta_{0}$, such that the trace of
${\mbox{\rv}}_{\mu}(t),$ $\mu\in(-\delta_{0},\delta_{0})$,
is a simple closed curve, let $\tilde{\mbox{\rv}}_{\mu}(t)$ be the map
${\mbox{\rv}}_{\mu}(t)$ considered  as a map with domain $S^{1}=\{ z\in\C; |\!|z|\!|=1\}$ (that is,
$\tilde{\mbox{\rv}}_{\mu}(e^{i\theta})={\mbox{\rv}}_{\mu}
(\frac{\tau_{\mu}}{2\pi}\theta), \,\tau_{\mu}=4t(\mbox{\vv}_{\mu},\mu)$).
Note that  $\tilde{\mbox{\rv}}_{\mu}$ and ${\mbox{\rv}}_{\mu}$ have the same trace.\\

From the continuous dependence of the  solutions  with respect to initial conditions,
and from the fact that  $t(\mbox{\vv},\mu)$
is continuous, it is straightforward to verify that, choosing $\delta_{0}$ and $\eta$
sufficiently small, we have  that, for $\,\mu\in(-\delta_{0},\delta_{0})$,
$\tilde{\mbox{\rv}}_{\mu}$ is close to $\tilde{\mbox{\rv}}_{0}$ in the space of 
maps from $S^{1}$ to $\R^{2}$ (with the $C^{1}$ topology).\\

Since $\tilde{\mbox{\rv}}_{0}$  is a embedding, and the space of embeddings is
open in the $C^{1}$ topology, we have that we can choose $\delta_{0}$ and $\eta$
sufficiently small  such  that, for $\,\mu\in(-\delta_{0},\delta_{0})$,
$\tilde{\mbox{\rv}}_{\mu}$ is 
also a embedding. This implies that the trace of $\tilde{\mbox{\rv}}_{\mu}$
is homeomorphic to $S^{1}$, that is, it is a simple closed curve in $\R^{2}.$\\

To show that we can choose $\delta_{0}$ sufficiently small such that  the trace of ${\mbox{\rv}}_{\mu}(t),$
$\,\mu\in(-\delta_{0},\delta_{0})$, encloses the origin, recall first that the trace of
${\mbox{\rv}}_{0}$ encloses the origin.
This means  that  $(0,0)$ is in the bounded component of
$\R^{2}-$(trace ${\mbox{\rv}}_{0}$)$=\R^{2}-C$.
Equivalently, ${\tilde{\mbox\rv}}_{0} :S^1
\rightarrow \R^2 -\{(0,0)\}$ is not homotopy trivial.\\

Choosing $\delta_{0}$ and $\eta$
sufficiently small we have that $\tilde{\mbox{\rv}}_{\mu}$ is close to $\tilde{\mbox{\rv}}_{0}$,
$\,\mu\in(-\delta_{0},\delta_{0})$.
Therefore ${\tilde{\mbox\rv}}_{\mu} :S^1 \rightarrow \R^2 -\{(0,0)\}$ is 
homotopic to ${\tilde{\mbox\rv}}_{0}$ in $\R^2 -\{(0,0)\}$. Hence
${\tilde{\mbox\rv}}_{\mu}$ is not homotopy
trivial either. It follows that the trace of ${\mbox\rv}_{\mu}$ encloses the
origin, $\mu\in(-\delta_{0},\delta_{0})$.
 \CaixaPreta

\vspace{0,8cm}

For the proof of the addendum we will use the following lemma. 
We  present a proof of this lemma in an appendix.

{\lem Let $[a,b],\,[c,d]\subset\R$ be closed intervals and
${f}:[a,b]\times[c,d]\rightarrow\R$ be continuous, such  that 
$$\left\{\begin{array}{l}
f(a,y)<0,\,\,\,y\in [c,d],\\ \\
f(b,y)>0,\,\,\,y\in [c,d].
\end{array}\right.$$

Then there is ${\cal W}\subset f^{-1}(0)$  connected and compact such that
${\cal W}_{y}:={\cal W}\cap \bigl( [a,b]\times\{y\} \bigr)\neq\emptyset,$ for all
$y\in [c,d].$
\label{3.1.2}}

\vspace{0,25cm}

\noindent {\bf Proof of the Addendum:} We use the notation from the proof of the theorem.

Recall that the  solutions  ${\mbox{\rv}}_{\mu}(t)=
{\mbox{\rv}}(t,{\mbox{\xx}}_{0},{\mbox{\vv}}_{\mu},\mu)$
are such that $l({{\mbox{\vv}}}_{\mu},\mu)=0$.
Moreover, if $l({{\mbox{\vv}}},\mu)=0$,
then  $({{\mbox{\vv}}},\mu)$ determines a periodic solution 
${\mbox{\rv}} _{ {} _{{\mbox{\vv}}},{}_{\mbox{$\mu$}}}$ of
$\stackrel{..}{\mbox{{\mbox{\rv}}}} \,= g({\mbox{{\mbox{\rv}}}},\mu )$.\\

Recall also that, by the claim above, there is ${\mbox{\vv}}_{+},\,{\mbox{\vv}}_{-}$,
$\delta_{0}>0,$ such that 
$l({\mbox{\vv}}_{+},\mu)>0$, $l({\mbox{\vv}}_{-},\mu)<0$, for
$\mu\in(-\delta_{0},\delta_{0})$.\\

Since ${\mbox{\vv}}_{+}\in V_{\delta}$, we can choose $\epsilon_{+}>0$ such that 
${\mbox{\vv}}_{+}= (1+\epsilon_{+})\,{\mbox{\vv}}_{0}$ 
(because $|\!|{\mbox{\vv}}_{0}|\!|<|\!|{\mbox{\vv}}_{+}|\!|$ and ${\mbox{\vv}}_{+}$
and ${\mbox{\vv}}_{0}$ have the same direction).
Analogously, we can choose $\epsilon_{-}<0$ such that ${\mbox{\vv}}_{-}= (1+\epsilon_{-})\,{\mbox{\vv}}_{0}.$ \\

Define $f: [\epsilon_{-},\epsilon_{+}]\times[0,\delta_{0}/2] \rightarrow\R,$
$f(s,\mu )=l((1+s){\mbox{\vv}}_{0},\mu)$.
Since $f(\epsilon_{+},\mu)>0,\,f(\epsilon_{-},\mu)<0,$
$\mu\in[0,\delta_{0}/2],$ the lemma above tells us that there is a connected compact 
${\cal W}\subset f^{-1}(0),$ such  that $${\cal W}_{\mu}:= {\cal W}\cap \Bigl(
 [\epsilon_{-},\epsilon_{+}]\times \{\mu\}\Bigr) \neq\emptyset$$

\noindent for $\mu \in[0,\delta_{0}/2]$. Therefore, for
$(\mbox{\vv},\mu)\in {\cal V}:=\Bigl\{ \,\Bigl(
(1+s){\mbox{\vv}}_{0},\mu\Bigr)\,;\, (s,\mu)\in {\cal
W}\,\Bigr\},$ we have that  $l(\mbox{\vv},\mu)=0$ and $\cal V$ satisfies
the conditions of the addendum to the theorem. \CaixaPreta

\vspace{0,2cm}

{\obs {\rm  For the proof of the existence of figure eight periodic solutions in \cite{AO}
we needed that the  solutions 
${\mbox{\rv}} _{ {} _{{\mbox{\vv}}_{\mu}},{}_{\mbox{$\mu$}}}$ given in theorem
\ref{3.1.1} above satisfy the following
property: ${\mbox{\rv}} _{ {}
_{{\mbox{\vv}}_{\mu}},{}_{\mbox{$\mu$}}}$
 intersects the $x$-axis in exactly two points. (It follows that these points
are ${\mbox\xx}_0$ and $-{\mbox\xx}_0$.)
To verify that the solutions  ${\mbox{\rv}} _{ {}
_{{\mbox{\vv}}_{\mu}},{}_{\mbox{$\mu$}}}$
satisfy this property, it is enough to apply proposition  \ref{3.0.17}
to the solutions  $\tilde{\mbox{\rv}}_{\mu}(e^{i\theta})={\mbox{\rv}}_{\mu}
(\frac{\tau_{\mu}}{2\pi}\theta)$ twice (see the end of the proof of theorem \ref{3.1.1}): 
once to $\{\, e^{i\theta}\,;\,-\frac{\pi}{2} \leq \theta \leq \frac{\pi}{2}\,\}$ and then to
$\{\, e^{i\theta}\,;\,\frac{\pi}{2} \leq \theta \leq \frac{3\pi}{2}\,\}$.}}

\vspace{0,4cm}

\section{Other Perturbations.}
 In this section we prove the existence of periodic solutions of perturbed 
problems. These solutions are close to a given circular solution  of the unperturbed problem
$\stackrel{..} {\mbox{\rv}} \,=-\frac{\mbox{$\kappa$}}{|\!|\mbox{{\mbox{\rv}}} |\!|^{\alpha
+2} }{\mbox{\rv}}\, ,$  with  $0\leq \alpha $, $\alpha\neq 1$ and $\kappa > 0$. We will require the
perturbed problems to be symmetric with respect to the $x$-axis. \\

Before proving the main result of this section let us recall some general facts about central 
force problems (see \cite{G}, p. 70-81, \cite{L}, p. 30-35, or \cite{Arn},
p. 33-41). \\

Let  $U:(0,+\infty)\rightarrow \R,$ $U\in C^{\infty}$. Consider the planar problem\\

\begin{equation}\stackrel{..}{\mbox{\rv}} \,=-\nabla
U(\mbox{\rv}) \label{UU}
\end{equation}\\

\noindent where  ${\mbox{\rv}}=(x,y)$ and
$U(\mbox{\rv})=U(r),\,r=|\!|\mbox{\rv}|\!|=\sqrt{x^{2}+y^{2}}.$
\vspace{0,2cm}
Using polar coordinates $(r,\varphi)$, the
Hamiltonian and the angular momentum can be written in the following form
\vspace{0,15cm}
\begin{equation}
\left\{
\begin{array}{lll}
H(r,\dot r)&=&\frac{1}{2}{\dot r}^{2} +\frac{K^{2}}{2r^{2}}+U(r),
\label{6.2}
\\&&\\
\,\,\,K &=& r^{2}\dot\varphi.
 \end{array}\right.
\end{equation}

\vspace{0,15cm}

Let $r(t)$ be a solution of this problem.\\

\noindent {\bf Remarks.}

\noindent (1) The following claims are equivalent: $\left\{
\begin{array}{lll}
i)\, \,\,\dot r(t_{0}) &=& 0, \\&&\\
ii)\, \,\,\mbox{\rv}(t_{0}) & \bot &\dot{\mbox{\rv}}(t_{0}).
\end{array}\right.$

\vspace{0,5cm}

\noindent (2) Suppose that there are $t_{1}<t_{2}$ such that $\dot r(t_{1})=\dot
r(t_{2})=0$, $ r(t_{1}) \neq r(t_{2})$ and 
 $\dot r(t)\neq 0,\,\,t\in (t_{1},t_{2}).$ Then we have two possibilities:

 a) $\dot r(t)>0,$ for all   $t\in (t_{1},t_{2}).\,\,$ It follows that 
$\left\{\begin{array}{lll}
 r(t_{1})= r_{min}
\\&&\\
 r(t_{2})= r_{max}
\end{array}\right.$
\vspace{0,2cm}

 b) $\dot r(t)<0,$ for all   $t\in (t_{1},t_{2}).\,\,$ It follows that 
$\left\{\begin{array}{lll}
 r(t_{1})= r_{max}
\\&&\\
 r(t_{2})= r_{min}
\end{array}\right.$
\vspace{0,2cm}

Certainly, if $\dot r\equiv 0$  then  $r_{max}=r_{min}=r(t)$
for all   $t$, and the solution $\mbox{\rv}(t)$ is circular.
The points where  $r=r_{min}$  are  called {\bf pericenters}, and the points where  $r=r_{max}$  are  
called {\bf apocenters}.\\

Therefore, if the orbit $\mbox{\rv}(t)$ is not circular and has, at least, one pericenter and one 
apocenter then the orbit describes a curve that goes from an apocenter to a pericenter or from a 
pericenter to an apocenter successively (see figure below). Moreover, $\mbox{\rv}(t)$
is defined for all $t$ and the movement happens in the interior of an annulus defined by the circles 
with radii $r_{min}$ and $r_{max}$. Consequently, if we know $r(t)$ between a
pericenter and an apocenter (o vice-versa), we know the whole 
function $r(t)$ (see fig. 3.8). The proof of these facts uses the symmetry of 
the problem and is similar to the proofs
of lemmas \ref{3.0.14} and \ref{3.0.15} (see \cite{Arn}, \cite{G}).

\vspace{0,4cm}

\begin{figure}[!htb]
\centering
\includegraphics[scale=0.5]{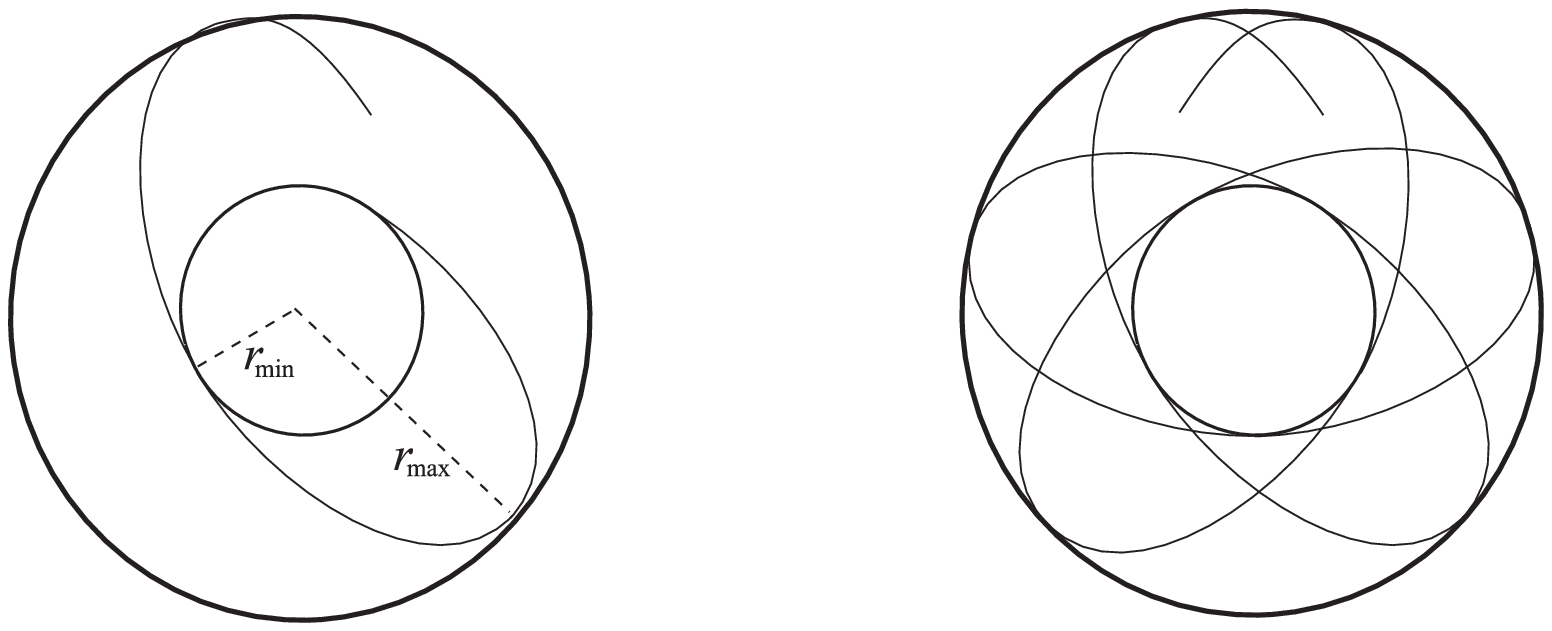}
\caption{\scriptsize{Orbit of $\mbox{\rv}(t)$}}
\end{figure}

\vspace{0,4cm}

\begin{figure}[!htb]
 \centering
 \includegraphics[scale=0.45]{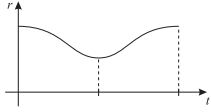}
 \caption{\scriptsize{Periodic function $r$}}
 \end{figure}

Note that if $\mbox{\rv}(t)$ has an apocenter and a pericenter then  $|\dot\varphi (t)|=
\frac{|K|}{\|\dot r (t)\|}\geq\frac{|K|}{r_{max}}>0$. Hence, since  $\mbox{\rv}(t)$ is defined
for all $t$, $lim_{t\rightarrow +\infty}|\varphi(t)|=+\infty,\,$ that is, $\,\mbox{\rv}(t)\,$
``goes around" the origin infinitely many times.\\

\noindent (3) The angle between a pericenter and a successive apocenter 
is given by
\vspace{0,15cm}
$$\Phi
=\int_{r_{min}}^{r_{max}}\frac{K/r^{2}\,dr}{\sqrt{2(E-U_{ff}(r))}}, $$

\noindent where  $U_{ff}=\displaystyle\frac{K^{2}}{2r^{2}}+U(r)$
is the effective potential and $ K$ is the angular momentum.\\

By the symmetry of the potential, we see that $\Phi$ does not
depend on which successive $r_{max},r_{min}$ we choose.\\

Let  $\{ \mbox{\rv}_{n}(t)\}$  be a sequence of solutions  that  approaches
a circular solution $ \mbox{\rv}_{0}(t)$, with radius
$r_{0}$. Suppose that  $ \mbox{\rv}_{n}(t)$ has an apocenter and a pericenter. Then  $\Phi _{n}\rightarrow
\pi\sqrt{\frac{U'(r_{0})}{3U'(r_{0})+r_{0}U''(r_{0})}},\,$ where 
$\Phi_{n}$  is the angle between a pericenter and a successive apocenter
of ${\mbox{\rv}}_{n}$ (see \cite{Arn}, p.37). Here ``$\{{\mbox{\rv}}_{n}(t)\}$ approaches 
$\mbox{\rv}_{0}(t)$" means that 
$( {\mbox{\rv}}_{n}(0),\dot {\mbox{\rv}}_{n}(0))\rightarrow
({\mbox{\rv}}_{0}(0),\dot{\mbox{\rv}}_{0}(0))$.\\

\noindent (4) Let $\mbox{\rv}_{0}(t)$ be a circular solution of (\ref{UU})
with ${\mbox{\rv}}_{0}(0)=\mbox{\xx}_{0}\neq 0\,$ and ${\dot{\mbox{\rv}}}_{0}(0)=\mbox{\vv}_{0}\neq 0$.
Assume $\mbox{\xx}_{0}\perp \mbox{\vv}_{0}$.
For $\epsilon >0$ let $\mbox{\rv}_{\epsilon}(t)$ denote a solution of (\ref{UU})
with ${\mbox{\rv}}_{\epsilon}(0)=\mbox{\xx}_{0}\,$ and ${\dot{\mbox{\rv}}}_{\epsilon}(0)=
(1+\epsilon)\mbox{\vv}_{0}$. Write $a=|\!|  \mbox{\xx}_{0} |\!|$ and assume that $U'(a)>0$.\\

\noindent{\bf Claim:} {\it 
If $\epsilon >0$ we have that $\mbox{\xx}_{0}$ is a
pericenter of $\mbox{\rv}_{\epsilon}(t)$; similarly, if $\epsilon <0$ we have that $\mbox{\xx}_{0}$
is an apocenter of $\mbox{\rv}_{\epsilon}(t)$.}\\

\vspace{0,15cm}

\begin{figure}[!htb]
\centering
\includegraphics[scale=0.45]{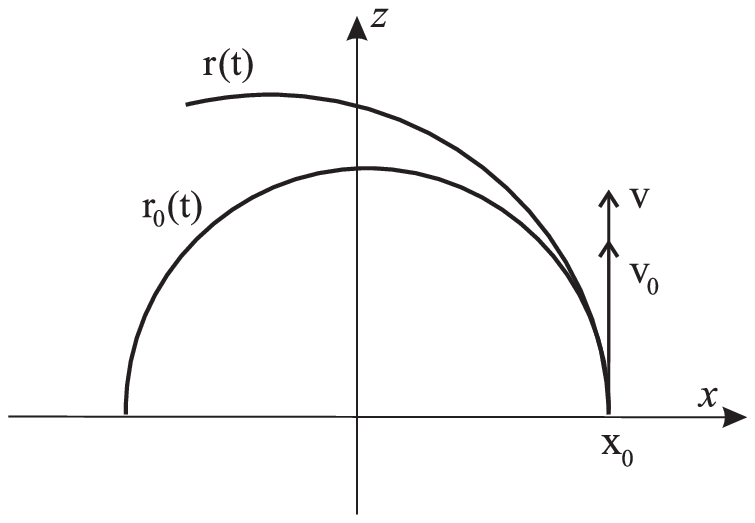}
 \caption{\scriptsize{}}
\end{figure}

Write $r_\epsilon=|\!|\mbox{\rv}_{\epsilon}|\!|$
and $v_{\epsilon}=|\!|\mbox{\vv}_{\epsilon}|\!|$, $\epsilon\in\R$. Then $v_{\epsilon}=(1+\epsilon)v_0$.
Write also $K_{\epsilon}$ for the angular momentum of $\mbox{\rv}_{\epsilon}$.
A simple calculation shows that $K_{\epsilon}=(1+\epsilon)av_0$.
Since $\mbox{\rv}_{0}$ is a circular solution we have $v_0^2=aU'(a)$.
Thus $K_{\epsilon}^2=(1+\epsilon)^2a^3U'(a)$.\\

Since $r_{\epsilon}$ is a solution of (\ref{6.2}) we have that 
$r_{\epsilon}$ satisfies (just differentiate the first equation):
$\stackrel{..}{r} \,=\frac{K^2}{r^3}-U'(r)$. Then, setting $t=0$ we have: 

$$\stackrel{..}{r}_\epsilon (0)\, =\,\frac{K_\epsilon^2}{a^3}-U'(a)\, =\, \frac{(1+\epsilon)^2a^3U'(a)}{a^3}
-U'(a)\, =\, \epsilon(2+\epsilon)U'(a)$$

It follows that the sign of $\stackrel{..}{r}_{\epsilon} (0)$ is equal to the sign of
$\epsilon$. Therefore, since $r_{\epsilon}(0)=a$, $\dot{r}_{\epsilon}(0)=0$ we have that,
if $\epsilon >0$,  $r_{\epsilon}(t)>a$ for $t\neq 0$ close to $0$; hence
$\mbox{\xx}_{0}$ is a pericenter of $\mbox{\rv}_{\epsilon}(t)$. Similarly for $\epsilon <0$.
This proves the claim.\\

\noindent (5) If we consider (\ref{UU}) of the form \vspace{0,2cm}
$$\stackrel{..}{\mbox{{\mbox{\rv}}}}
\,=\,-\frac{\mbox{$\kappa$}}{|\!|\mbox{\rv}|\!|^{\alpha +2} }
{\mbox{\rv}},\,\,\,\mbox{ where }\,\,{\mbox{\rv}}=(x,y),\,\,\kappa >
0,\,\,\mbox{  with }\,\,\, 0\leq \alpha , $$

\noindent the potential $U$ of the problem is given by:
\vspace{0,2cm}
\begin{equation} U(r)\,=\, \left\{
\begin{array}{lll}
 -\displaystyle\frac{\gamma}{r^{\alpha}},\,\,\,\,\mbox{with}\,\,\,\,\alpha
\gamma=\kappa, &&\,\,0<\alpha ,
\\&&\\
 \kappa\,ln\,r, && \,\,\alpha =0.
\end{array}\right.
\label{U}
\end{equation}

Observe that, if ${\mbox{\rv}}_{0}(t)$  is a circular solution  of an attractive central problem,
with  potential $U$, we have that, for  ${\mbox{\rv}}_{0}(0)=\mbox{\xx}_{0}\,$ and
${\dot{\mbox{\rv}}}_{0}(0)=\mbox{\vv}_{0}:$
$$|\!|\mbox{\vv}_{0}|\!|=\sqrt{pU'(p)},\,\,\,\mbox{with}\,\,
p=|\!|\mbox{\xx}_{0}|\!|.$$ \vspace{0,2cm}
For $U$  as in (\ref{U}), if $0<\alpha ,$ $\,|\!|\mbox{\vv}_{0}|\!|=
\frac{\sqrt{\alpha \gamma}}{\sqrt {p^{\alpha}}}\,$ and if $\alpha
=0,$ $\,|\!|\mbox{\vv}_{0}|\!|=\sqrt{\kappa}$. Note that, in this last case,
$\,|\!|\mbox{\vv}_{0}|\!|$  is independent of $\mbox{\xx}_{0}$.\\

\vspace{0,6cm}

\noindent {\bf Theorem 0.2} {\it 
 Let  $C$ be a circle centered at the origin  $(0,0)$ of $\R^{2}$,
${\mbox{\xx}}_{0}\in C\cap \{ x-axis\}$ and $U$ an
open neighborhood of $C$ of the form $C\subset U\subset (\R^{2}-\{ (0,0)\})$.\\
Let $a>0$ and $g : U\times (-a,a)\rightarrow\R^{2}$ continuous such that\\
(i) $g(\mbox{{\mbox{\rv}}},0)=\,-\frac{\mbox{$\kappa$}}
{|\!|\mbox{\rv}|\!|^{\alpha +2}}
{\mbox{\rv}},$ where
$\mbox{{\mbox{\rv}}}=(x,y),\,\,\kappa > 0,$ and $\,\,0\leq \alpha 
,\,\,\, \alpha\neq 1$,\\
(ii) $g({\mbox{\rv}},\mu )$ is $C^{1}$ in the variable ${\mbox{\rv}}\in
U$, for each $\mu ,$\\
(iii) for all $\mu $, $g$ is invariant (as a vector field)
by the reflection
$$\varphi(x,y)=(x,-y).$$
\quad  Then there is $\delta_{0}$, $0<\delta_{0}<a$, with the 
following property.  For each $\mu\in
(-\delta_{0},\delta_{0})$ there is a velocity ${\mbox{\vv}}_{\mu}$ 
such that the solution
${\mbox{\rv}}_{{}_{{\mathrm{\mbox{\vv}}}_{\mu}},{}_{{\mbox{$\mu$}}}}(t)
:= {\mbox{\rv}} (t,{\mbox{\xx}}_{0},{\mbox{\vv}}_{\mu},\mu)$
of $\stackrel{..}{\mbox{{\mbox{\rv}}}} \,=g(\mbox{{\mbox{\rv}}},\mu)$ is periodic.
Moreover, given $\eta >0$ we can choose  $\delta_{0}>0$ such that

(1) the traces of these solutions
are simple closed curves, symmetric with respect to the x-axis, 
and enclose the origin,

(2)  all velocities ${\mbox{\vv}}_{\mu}$ are $\eta$-close, $\mu < 
\delta_{0}$.}

\vspace{0,4cm}

The Remarks 2.1 also hold in this case. Thus we also give an addendum to
theorem 0.2:\\

\vspace{0,2cm}

\noindent {\bf Addendum (to theorem \ref{3.2.1})}
{\it
We can choose  $\delta_{0}>0$ in theorem   \ref{3.2.1} such 
that there is a compact connected set ${\cal V}\subset \R^{2}\times\R$ 
 with the following properties:

\noindent 1) ${\cal V}_{\mu}:={\cal V}\cap \left( \R^{2}\times\{\mu\}\right)\neq
\emptyset,$ for all $\mu\in[0,\delta_{0}]$.

\noindent 2) ${\mbox{\rv}}_{{}_{\mbox{\vv}},{}_{\mbox{$\mu$}}}$ is 
periodic, for $(\mbox{\vv},\mu)\in {\cal V}.$

Moreover, the trace of ${\mbox{\rv}}_{{}_{\mbox{\vv}},{}_{\mbox{$\mu$}}}$ is a simple closed curve 
symmetric with respect to the x-axis, and encloses the origin.}\\

The proof of  theorem \ref{3.2.1} is similar to the proof of theorem 0.1 (Kepler's perturbation problem). 
We just need to study the behavior of the solutions near a circular solution of the unperturbed
problem (which we knew in the case of Kepler's problem: they are all ellipses). \\

Let  ${\mbox{\rv}}_{0}(t)=(x_0(t),y_0(t))$ be a circular solution of the planar problem 
\vspace{0,15cm}
\begin{equation}
\stackrel{..}{\mbox{{\mbox{\rv}}}} \,
=-\frac{\kappa}{|\!|\mbox{\rv}|\!|^{\alpha
+2}}{\mbox{\rv}},\,\,\,0\leq\alpha ,\,\,\,\alpha\neq 1\,\,\,\,\kappa>0
\label{5.4}
\end{equation}

Suppose that ${\mbox{{\mbox{\rv}}}}_{0}(0)={\mbox{\xx}}_{0}$ and
$\dot{\mbox{{\mbox{\rv}}}}_{0}(0)={\mbox{\vv}}_{0},
\,{\mbox{\vv}}_{0}$ with the same direction as the positive $y$-axis and ${\mbox{\xx}}_{0}$ in the 
positive $x$-axis.\\

Let  ${\mbox{\rv}}_{\epsilon}(t)=(x_{\epsilon}(t),y_{\epsilon}(t))$ be a solution of 
(\ref{5.4}) with  
the same initial position ${\mbox{\xx}}_{0}$ and with 
initial velocity ${\mbox{\vv}}_{\epsilon}=(1+\epsilon){\mbox{\vv}}_0$.\\

Write $r_{\epsilon}=|\!|{\mbox{\rv}}_{\epsilon}|\!|$,
$v_{\epsilon}=|\!|{\mbox{\vv}}_{\epsilon}|\!|$ and $v_{0}=|\!|{\mbox{\vv}}_{0}|\!|.$
Let $\varphi_{\epsilon}$ be such that $(r_{\epsilon},\varphi_{\epsilon})$ are 
the polar coordinates of
${\mbox{\rv}}_{\epsilon}$, with $\varphi_{\epsilon}(0)=0$.\\
 
Since $\dot{\mbox{\rv}}_{\epsilon}(0)$ is orthogonal to ${\mbox{\rv}}_{\epsilon}(0)$, 
we have $\dot r _{\epsilon}(0)=0$ 
(see remark (1) above) and $r_{\epsilon}(0)$  is a maximum or minimum of $r_{\epsilon}(t)$;
hence ${\mbox{\rv}}_{\epsilon}(0)$ is an apocenter or pericenter of the solution 
${\mbox{\rv}}_{\epsilon}(t).$ Note that remark (4) above implies that 
${\mbox{{\mbox{\rv}}}}_{\epsilon}(0)={\mbox{\xx}}_{0}$
is a pericenter of ${\mbox{\rv}}_{\epsilon}$, for $\epsilon >0$
and an apocenter of ${\mbox{\rv}}_{\epsilon}$, for $\epsilon <0$.\\

It follows from proposition \ref{3.0.18} that, for $\epsilon$ close to $0$, there is a
minimum $t_{\epsilon}^* >0$ such that ${\mbox{\rv}}_{\epsilon}$ intersects 
the negative $x$-axis (transversally) in ${\mbox{\rv}}_{\epsilon}( t_{\epsilon}^*)$.
Equivalently, $\varphi_{\epsilon}(t^*_{\epsilon})=\pi$.
The next lemma tells us when ${\mbox{\vv}}_{\epsilon}( t_{\epsilon}^*)$ points to the
left or to the right; i.e. when the sign of $\dot x _{\epsilon}( t_{\epsilon}^*)$
is positive or negative.

{\lem The sign of $\dot x _{\epsilon}( t_{\epsilon}^*)$ is given by the following table:
\label{3.2.2}}

\newpage

\begin{table}[h]
\centering
\begin{tabular}{|c||c|c|}
\hline
 & $0\leq\alpha<1$ & $1<\alpha$ \\
\hline \hline $ \epsilon>0$ & $\dot x _{\epsilon}( t_{\epsilon}^*)>0$ &
$\dot x _{\epsilon}( t_{\epsilon}^*)<0$ \\
\hline
 $ \epsilon<0$ & $\dot x _{\epsilon}( t_{\epsilon}^*)<0$ &
$\dot x _{\epsilon}( t_{\epsilon}^*)>0$ \\
\hline

\end{tabular}
\caption{Sign of $\dot x _{\epsilon}( t_{\epsilon}^*)$ }
\end{table}

\vspace{0.5cm}

\noindent{\bf Proof:}
Let us assume first that $0\leq\alpha < 1$. In this case, it is
known that, for small $\epsilon$, ${\mbox{\rv}}_{\epsilon}$ has an apocenter
and a pericenter (see \cite{G}). \\

Let $\Phi( \epsilon )$  be the angle between a  pericenter and a successive
apocenter of the solution ${\mbox{\rv}}_{\epsilon}(t)$. Since $0\leq \alpha <1$, we have that 
$\frac{\pi}{2}<lim_{\epsilon\rightarrow 0}\Phi( \epsilon )\,=\,\frac{\pi}
{\sqrt{2-\alpha}}<\pi$ (see remark (3) above). This, together with remark (4) above, imply
that, for $\epsilon$ sufficiently close to 0, we have: (see fig. 3.10)\\

(1) if $\epsilon >0$, ${\mbox{\rv}}_{\epsilon}(t_{\epsilon}^*)$ is between
an apocenter and a pericenter (in that exact order, where we consider one
``first'' if it happens at an earlier time). If follows that 
$r_{\epsilon}$ is a decreasing function in this interval; hence
$\dot{r}_{\epsilon}(t^*_{\epsilon})<0$. Since $x_{\epsilon}=r_{\epsilon}cos\, \varphi_{\epsilon}$
and $\varphi_{\epsilon}(t^*_{\epsilon})=\pi$ we have that
$\dot{x}_{\epsilon}(t^*_{\epsilon})=\dot{r}_{\epsilon}(t^*_{\epsilon})\cos \pi >0.$\\

\begin{figure}[!htb]
\centering
\includegraphics[scale=0.55]{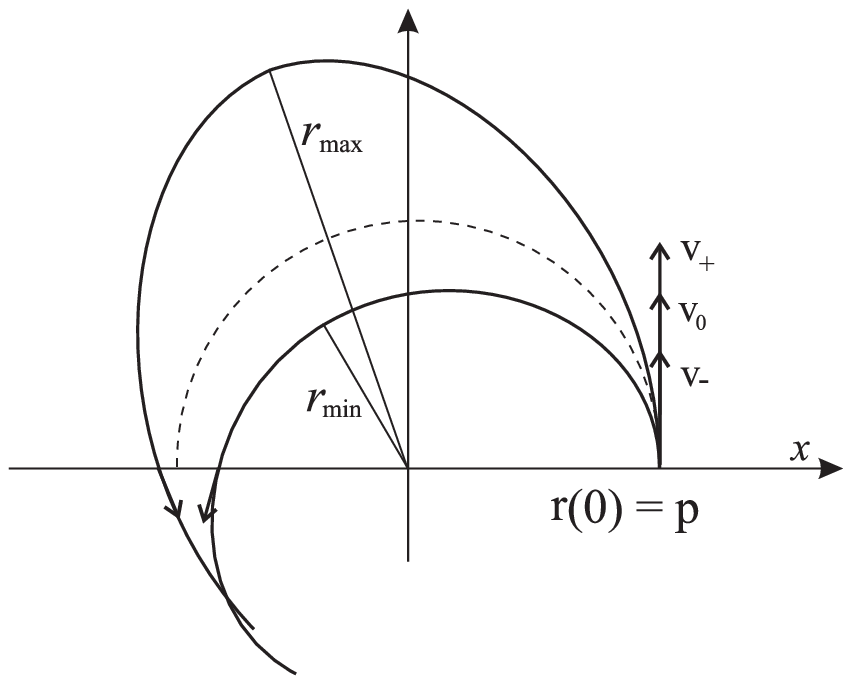}
\caption{\scriptsize{$0\leq\alpha<1$}}
\end{figure}

(2) If $\epsilon <0$, ${\mbox{\rv}}_{\epsilon}(t_{\epsilon}^*)$ is between
a pericenter and an apocenter. Then
$r_{\epsilon}$ is an increasing function in this interval; hence
$\dot{r}_{\epsilon}(t^*_{\epsilon})>0$. As before, it follows that
$\dot{x}_{\epsilon}(t^*_{\epsilon})<0.$\\

Now, let us assume that $1<\alpha$. By the form of the graph of the effective potential
in this case (see \cite{G}) we have that, for $t>0$,  $r_{\epsilon}(t)$ is either 
an increasing function or a decreasing function. By remark (4) above, $r_{\epsilon}(t)$
is an increasing function, for $\epsilon >0$ and $t>0$, and 
 $r_{\epsilon}(t)$ is a decreasing function, for $\epsilon <0$ and $t>0$.
An argument similar to the one used in cases (1) and (2) above completes the proof
of the lemma. \CaixaPreta\\

\noindent {\bf Proof of theorem \ref{3.2.1}:}
As before, ${\mbox{\rv}}(t,{\mbox{\xx}},{\mbox{\vv}},\mu)$ denotes
a solution of  $ \stackrel{..}{\mbox{{\mbox{\rv}}}}
\,= g({\mbox{{\mbox{\rv}}}},\mu ),$ with initial conditions
${\mbox{\rv}}(0,{\mbox{\xx}},{\mbox{\vv}},\mu)={\mbox{\xx}}$ and
$\dot{\mbox{\rv}}(0,{\mbox{\xx}},{\mbox{\vv}},\mu)={\mbox{\vv}}$.\\

We proceed in a similar fashion as in the final part of the proof
of theorem \ref{3.1.1}\\

Choose $\bar t >0$ such  that  ${\mbox{\rv}}_{0}(t)=
{\mbox{\rv}}(t,{\mbox{\xx}}_{0},{\mbox{\vv}}_{0},0)$,
restricted to the interval $[0,\bar t\,]$, intersects (transversally) the 
negative $x$-axis in a a single point.\\

Let   $\delta>0$, $\delta$
sufficiently small (as in proposition  \ref{3.0.18},  with  $\bar t$  as above and $E= \{(x,0);x\leq 0\}$), 
and define $V_{\delta} = \{\sigma {\mbox{\vv}}_{0};\,\sigma\in(1-\delta,
1+\delta)\}\subset\R^{2}$.
Define the function $\tilde{l}
:V_{\delta}\times (-\delta,\delta)\rightarrow
\R,\,\,  \tilde{l}({\mbox{\vv}},\mu)=\dot{x}(t({\mbox{\vv}},\mu),\mbox{\xx}_{0},
{\mbox{\vv}},\mu)$, where 
$\dot{x}(t({\mbox{\vv}},\mu),\mbox{\xx}_{0},{\mbox{\vv}},\mu)$
is the first coordinate of the velocity vector
$\dot{\mbox{{\mbox{\rv}}}}(t({\mbox{\vv}},\mu),\mbox{\xx}_{0},
{\mbox{\vv}},\mu)=(\dot{x}(t({\mbox{\vv}},\mu),
\mbox{\xx}_{0},{\mbox{\vv}},\mu),\dot{y}(t({\mbox{\vv}},\mu),
\mbox{\xx}_{0},{\mbox{\vv}},\mu))$ and $t({\mbox{\vv}},\mu)>0$
is the minimal time at which the solution intersects the negative $x$-axis.
Here $t({\mbox{\vv}},\mu):= t({\mbox\xx}_{0},{\mbox{\vv}},\mu)$  is as in
proposition  \ref{3.0.18}.\\

As in claim 2.1, we can verify, using lemma 3.1, that there are
$\delta_{0}$, $0<\delta_{0}<\delta$, and ${\mbox{\vv}}_{-},$
${\mbox{\vv}}_{+}\in V_{\delta}$, $|\!|{\mbox{\vv}}_{-}|\!|<|\!|{\mbox{\vv}}_{0}|\!|<|\!|{\mbox{\vv}}_{+}|\!|$,
such that, for all $\mu\in (-\delta_{0},\delta_{0}),$ we have: \\

(1) If $0\leq\alpha<1$, then
$\tilde{l}({{\mbox{{\mbox{\vv}}}}}_{-},\mu)<0$ and $\tilde{l}({{\mbox{{\mbox{\vv}}}}}_{+},\mu)>0$.

(2) If $1<\alpha$, then
$\tilde{l}({{\mbox{{\mbox{\vv}}}}}_{+},\mu)<0$ and $\tilde{l}({{\mbox{{\mbox{\vv}}}}}_{-},\mu)>0$.\\

In any case, by the intermediate value theorem, there is
$\delta_{0}>0$ such that, for all  $\mu\in(-\delta_{0},\delta_{0})$, there 
exists ${{\mbox{\vv}}}_{\mu}\in V_{\delta}$ such  that 
$\tilde{l}({{\mbox{\vv}}}_{\mu},\mu)=0$.
Therefore, for each  $\mu\in (-\delta_{0},\delta_{0}),$ the solution of the perturbed problem 
${\mbox{\rv}}_{\mu}(t):={\mbox{\rv}}_{{}_{{\mbox\vv}_{\mu}},{}_{\mbox{$\mu$}}}(t)=
{\mbox{\rv}}(t,{\mbox{\xx}}_{0},{\mbox{\vv}}_{\mu},\mu)$,
intersects orthogonally the $x$-axis at $t=t({{\mbox{\vv}}}_{\mu},\mu).$\\

Finally,  lemma \ref{3.0.14} implies that  ${\mbox{\rv}}_{\mu}(t)$ can be extended to a
periodic solution,  with  period $2t({{\mbox{\vv}}}_{\mu},\mu)$,
whose trace is symmetric  with respect  to the $x$-axis.
The rest of the proof is exactly the same  as the proof of theorem \ref{3.1.1}. \CaixaPreta

\vspace{1cm}

The proof of the addendum to theorem \ref{3.2.1}  is similar as the proof of the addendum
to theorem \ref{3.1.1}.
\vspace{0,4cm}

{\obs {\rm  For the proof of the existence of figure eight periodic solutions
in \cite{AO} we need also that the  solutions 
${\mbox{\rv}} _{ {} _{{\mbox{\vv}}},{}_{\mbox{$\mu$}}}$ given by
theorem \ref{3.2.1} above, satisfy also the following
property: ${\mbox{\rv}} _{ {}
_{{\mbox{\vv}}_{\mu}},{}_{\mbox{$\mu$}}}$
 intersects the $x$-axis in exactly two points. One of these 
 points is ${\mbox\xx}_0 =(x_0, 0)$ and the other is $(x^{'}_0,0)$, and we can take
$x_0 >0$ and $x'_0<0$.
To verify that the  solutions  ${\mbox{\rv}} _{ {}
_{{\mbox{\vv}}_{\mu}},{}_{\mbox{$\mu$}}}$
satisfy this property, it is enough to apply twice proposition  \ref{3.0.17}
to the solutions 
$\tilde{\mbox{\rv}}_{\mu}(e^{i\theta})={\mbox{\rv}}_{\mu}
(\frac{\tau_{\mu}}{2\pi}\theta)$. Once to $\{\,
e^{i\theta}\,;\,-\frac{\pi}{2} \leq \theta \leq \frac{\pi}{2}\,\}$ and then to
$\{\, e^{i\theta}\,;\,\frac{\pi}{2} \leq \theta \leq \frac{3\pi}{2}\,\}$.}\label{pp2}}

\vspace{0,4cm}

\appendix

\section{ Proof of Proposition \ref{3.0.17}:}

Consider $E$ contained in the $y$-axis.
Let  $\pi _{1}(x,y)=x$ and let  $\alpha(t_{0})$ be the unique point
in the intersection of $\alpha$  with  $E$. 
Note that $\alpha(t)\in E$ if and only if $\pi_{1}(\alpha (t))=0$.
We assume $ \pi_{1}({\dot \alpha}(t))>0$.\\

Since ${\dot \alpha}$  is continuous, there is an interval $[a,b]\subset
(0,\bar t\,),\,\,t_{0}\in (a,b),\,$ such  that  $\,\pi_{1}({\dot \alpha}(t))>0,$
for all   $t\in [a,b].$ Also, we can suppose that 
$diam(\alpha[a,b])<\frac{1}{3}\,dist(\alpha(t_{0}),\pat E)=\gamma.$
Hence the function $\pi_{1}\alpha (t)$  is an increasing function on $[a,b]$. Then $
\pi_{1}\alpha (a)<0\,$ and $\,\pi_{1}\alpha (b)>0$, because $\pi_{1}\alpha (t_0) =0.$\\

Let  $\varepsilon_{1}=min\{dist(\alpha(t),E);\,t\in [0,\bar
t\,]\setminus
(a,b)\}\,$, $\varepsilon_{2}=min\{\pi_{1}{\dot \alpha }
(t);\,t\in[a,b]\}$ and consider $\epsilon
=\frac{1}{2}min\{\varepsilon_{1},\varepsilon_{2},\gamma \}$. (Note that  
$\varepsilon_{1},\varepsilon_{2}$ and $\gamma$  are  positive.)\\

Let  $\beta :[0,\bar t\,]\rightarrow\R^{2}$ be of class $C^{1}$,
such  that  $\|\alpha -\beta\|<\epsilon,\,\,\|\dot\alpha
-\dot\beta\|<\epsilon$. We have the following claims:
\vspace{0,2cm}

\noindent {\bf Claim 1:} $\beta (t)\notin E,\,\,t\in[0,\bar
t\,]\setminus (a,b).$

In fact, for all   $x\in E$, we have  $dist(x,\beta)\geq
dist(x,\alpha)-dist(\alpha, \beta)\geq dist(E,\alpha)-dist(\alpha,
\beta),$ and since $t\in [0,\bar t\,]\setminus (a,b),$ we have
$dist(x,\beta)\geq
\varepsilon_{1}-\frac{\mbox{$\varepsilon_{1}$}}{2}=
\frac{\mbox{$\varepsilon_{1}$}}{2}.$  Since this holds for all  
$x\in E$, we have  $dist(E,\beta)\geq
\frac{\mbox{$\varepsilon_{1}$}}{2}>0.$ \vspace{0,2cm}

\noindent {\bf Claim 2:} There exists a single $t_{1}\in (a,b)$ such 
that  $\beta (t_{1})\in E.$ Moreover, $\beta(t_{1})$  is an interior point of
$E$ and ${\dot \beta}( t_{1})\notin E$.\\

First we verify that  $\pi_{1}\beta(t)$ is an increasing function on $[a,b]$
and that  $\pi_{1}\beta(a)<0$ and $\pi_{1}\beta(b)>0$.
In fact, since $|{\dot \alpha}(t)-{\dot \beta}(t)|<\frac{\mbox{$\varepsilon_{2}$}}{2}$,
we have that  $|\pi_{1}{\dot \alpha}(t)-\pi_{1}{\dot
\beta}(t)|<\frac{\mbox{$\varepsilon_{2}$}}{2}$,
then
\vspace{0,17cm}
$$\pi_{1}{\dot \beta}(t)>\pi_{1}{\dot
\alpha}(t)-\frac{\varepsilon_{2}}{2}>0,\,\,\,\mbox{for}
\,\,\,t\in [a,b].$$

Also,
$|\beta(t)-\alpha(t)|<\frac{\mbox{$\varepsilon_{1}$}}{2},$ then
$|\pi_{1}\beta(t)-\pi_{1}\alpha(t)|<\frac{\mbox{$\varepsilon_{1}$}}{2}$.
For $t=a$, we have:  \vspace{0,17cm}
$$\pi_{1}\beta(a)<\pi_{1}\alpha(a)+\frac{\varepsilon_{1}}{2}<
\pi_{1}\alpha(a)+\frac{|\pi_{1}\alpha(a)|}{2}<0.$$

\noindent In the same way, we verify that $\pi_{1}\beta(b)>0$.
We have then  that  there exists a single  $t_{1}\in (a,b)$  such that 
$\pi_{1}\beta(t_{1})=0$, that is, $\beta(t_{1})\in \{\, y\}$-axis, not
necessarily in $E$. On the other hand,
\vspace{0,17cm}
$$|\beta (t_{1})-\alpha(t_{0})|\leq |\beta
(t_{1})-\alpha(t_{1})|+ |\alpha (t_{1})-\alpha(t_{0})| \leq
\epsilon +\gamma \leq 2\gamma =\frac{2}{3}dist (\alpha(t_{0}),\pat
E).$$
\vspace{0,2cm} 
\noindent that is, there is a single 
$t_{1}\in (a,b)$  with  $\beta(t_{1})$ in the interior $int\, E$ of $E$.
Since $\pi_{1}{\dot \beta}(t_{1})\neq 0$, we have  that 
${\dot \beta}(t_{1})\notin E.$ This proves claim 2. \\

From claims 1 and 2, we conclude that  there is a single 
$t_{1}\in [0,\bar t\,]$  with  $\beta(t_{1})\in\,E$. Moreover,
$\beta(t_{1})\in\, int\, E$, ${\dot \beta}(t_{1})\notin
E$ and $t_{1}\in(0,\bar t)$. This proves the proposition.
\CaixaPreta

\section{ Proof of lemma  \ref{3.1.2}:}
\quad Here we present a proof of lema \ref{3.1.2}. 
The proof is a direct application of the naturality of the
Mayer-Vietoris sequence. The referee pointed out to us
that a proof using Leray-Schauder degree theory can be found in \cite{R}.

\vspace{0,4cm}

\noindent {\bf Proof of the lemma 2.2:} Let $Z=[a,b]\times [c,d]$, $X=f^{-1}(0)$,
$A= \{\, a\}\times [c,d],$  $B=\{\,b\,\} \times [c,d]$, $Y_0 =[a,b]\times\{\, c\, \}$,
 $Y_1 =[a,b]\times \{\,d\,\}$.
 We want to prove that there is a compact connected ${\cal W}\subset X,$ with
 ${\cal W}\cap Y_i \neq \emptyset,$ $i=0,1$. If this is not the case, a result 
in general topology asserts that there are disjoint compact sets $X_0, X_1 \subset X$ with
 $Y_i \cap X\subset X_i$ and $X=X_0 \cup X_1$.
Note that $ H_{0}(A\cup B) \cong \Z_{2} \oplus \Z_{2}$ ($H_* $ denotes singular 
homology with $\Z_2$ coefficients). We identify $(1,0) \in \Z_2 \oplus \Z_2$ with the 
class in 
$H_{0}(A\cup B)$ determined by $A$, and $(0,1)$ with the class determined by $B$.
Write $x=(1,1)\in H_0(A\cup B)$. Since $A, B$ belong to the same path-connected component of 
$Z- X_i$ (because $Y_{i} \subset Z- X_{i+1},\,\,i\, (mod\, 2)$),
we have that $i_i (x) =0\in H_0 (Z- X_i),$ where $i_i:H_{0}(A\cup B)
\rightarrow H_{0}(Z- X_i)$ is induced by the inclusion.
Note that, since $f(A)\subset (-\infty,0),\,f(B)\subset (0,+\infty),$ $i(x)\neq 0\in 
H_{0}(Z- X)$, where $i:H_{0}(A\cup B) \rightarrow H_{0}(Z- X)$ 
is induced by the inclusion.
Consider the following diagram of Mayer-Vietoris sequences:\\

 \xymatrix{
\ar[r] &H_{1}(A\cup B)\ar[r]\ar[d] & H_{0}(A\cup B)
\ar[d]^{\iota}\ar[r]^{\phi\,\,\,\,\,\,\,\,\,\,\,\,\,\,\,\,\,} &
H_{0}(A\cup B) \oplus H_{0}(A\cup B)
 \ar[d]^{\iota_{1}\oplus \iota_{2}}\ar[r]  &H_{0}(A\cup B)\ar[r]\ar[d] &0\\
\ar[r] &H_{1}(Z)\ar[r]^{\theta\,\,\,\,\,\,\,\,\,\,} &
H_{0}(Z-X)\ar[r]^{\psi\,\,\,\,\,\,\,\,\,\,\,\,\,\,\,\,\,\,\,\,\,\,\,\,\,}
&
 H_{0}(Z-X_{0}) \oplus H_{0}(Z-X_{1})
\ar[r]  &H_{0}(Z)\ar[r] &0\\}
\vspace{.2in}

The first line is the Mayer-Vietoris sequence of  $A\cup B=(A\cup B)\cup
(A\cup B)$, and the second line is the Mayer-Vietoris sequence of
$Z=\,(Z-X_{0})\,\cup\,(Z-X_{1}).\,$ Note that
$(Z-X_{0})\,\cap\,(Z-X_{1})=Z-(X_{0}\cup X_{1})= Z-X$.
The vertical maps are induced by inclusions. 
Since  $\phi(x)=(x,x)$, we have
$(\iota_{1}\oplus \iota_{2})\,\phi (x)= (\iota_{1}(x),
 \iota_{2}(x))= (0,0)$.
Then $\psi\,\iota ({x})=0$ and we have $\iota (x)\in Im\,\theta.$
But $H_{1}(Z)=0.$
Consequently, $\iota (x)=0,$
a contradiction.   \CaixaPreta

{}

\end{document}